\documentclass[12pt,a4paper]{article}
\usepackage{latexsym}
\usepackage{amsfonts}
\usepackage{amssymb}
\usepackage{calrsfs}
\usepackage{amsmath}
\newtheorem{theorem}{Theorem}[section]
\newtheorem{lemma}{Lemma}[section]

\newtheorem{proposition}{Proposition}[section]

\newtheorem{definition}{Definition}[section]

\newtheorem{remark}{Remark}[section]

\numberwithin{equation}{section}

\usepackage{rotating}
\usepackage{wrapfig}
\usepackage{graphicx}

\newcommand{\R}{\mathbb R}

\newcommand{\N}{\mathbb N}

\title{Hyperbolic Volterra equations of convolution type in Sobolev spaces}

\author{Nadezhda~A. Rautian and Victor~V. Vlasov}
\date{}

\begin{document}
\maketitle

\begin{abstract}
We study the correct solvability of an abstract integro-differential equations in Hilbert space generalizing integro-differential equations arising in the theory of viscoelastisity. The equations under considerations are the abstract hyperbolic equations perturbed by the terms containing Volterra integral operators. We establish the correct solvability in the weighted Sobolev spaces of vector-valued functions on the positive semiaxis.
\end{abstract}

\footnotetext {The authors were financially supported by Russian Foundation for Basic Research (grant N 14-01-00349), and also Russian Scientific Foundation (grant N 14-11-00754).}

\section{Introduction}

Numerous problems arising in applications come the to research of the integro-differential equations. Let's point to some problems of mechanics and the physics which studying is reduced to the integro-differential equations considered in present paper.

The first class of problems are the problems arising in the theory of viscoelasticity. Integral terms like convolution describe in this case long-term memory, thus functions of a kernel of convolution are defined as a result of experiment. Curves often received as a result of experiment are approximated in practice by the sum of finite number of exponentials or series of exponents:
$$
K(t) \sim \sum\limits_{k = 1}^{\infty}{ {c_k} e^{ - \gamma _k t} }. 
$$
Rather complete description of the problems arising in the theory of viscoelasticity is given in the recent monograph \cite{AFG} (see also  \cite{24},   \cite{37} and references therein).

The second class of problems are the problems of acoustics of emulsions. The mix of two liquids with various characteristics (density, viscosity, compressibility coefficient) is called an emulsion. It is possible to prove by the methods of the homogenezation theory that the equation describing an average value of sound pressure for one-dimensional distributing of a sound wave has the abstract form coinciding with a form of the equation, considered in this work (see  \cite{30},  \cite{45},  \cite{411}). 

The third class of problems are the problems of homogenezation in multiphase media where one of phases is elastic (or viscoelastic) media, and another is viscous (compressible or incompressible) liquid (see \cite{37}, \cite{30}).

This paper is devoted to researching of the integro-differential equations with unbounded operator coefficients in  Hilbert space. Most of the
equations under consideration are abstract hyperbolic equations perturbed by terms containing Volterra integral
operators. These equations are the abstract forms of the integro-differential equations arising in the theory of viscoelastisity (see \cite{24}, \cite{AFG}) and the Gurtin-Pipkin integro-differential equations (see \cite{20}, \cite{42}, \cite{6} for more details), which describes heat propagation in media; it also arises in homogenization problems in
porous media (Darcy law) (see \cite{45}, \cite{30}, \cite{35}).

It is shown that the initial boundary value problems for these equations are well-solvable in Sobolev spaces on the
positive semiaxis.

Due to the fact that we study not concrete particular partial integro-differential equation but a wide class of integro-differential equations, it is
natural and convenient to consider integro-differential equations with unbounded operator coefficients (abstract
integro-differential equations), which can be realized as integro-differential partial differential equations with
respect to spatial variables when necessary. For the self-adjoint positive operator~$A$ considered in what follows we
can take, in particular, the operator $A^2y=-y''$, where $x\in(0,\pi)$, $y(0)=y(\pi)=0$, or the operator $A^2y=-\Delta
y$ satisfying the Dirichlet conditions on a bounded domain with sufficiently smooth boundary or more general elliptic self-adjoint operator of the order $2m$ on a bounded domain with sufficiently smooth boundary. At present, there is an
extensive literature on abstract integro-differential equations (see, e.g., \cite{1}-\cite{3}, \cite{4}, \cite{5}, \cite{14}--\cite{21}, \cite{11}--\cite{18}, \cite{8}, \cite{AFG} and the references therein).

In \cite{1}--\cite{3}, \cite{4}, \cite{5}, \cite{10}--\cite{12}, \cite{8}, \cite{9} (see also the references therein), integro-differential equations with
principal part being an abstract parabolic equation were studied. Equations with principal part being an abstract
hyperbolic equation are significantly less studied. The works most closely related to this questions are
\cite{22}, \cite{14}, \cite{41}--\cite{18}.

We emphasize that our method of the proof of the theorem on the correct solvability of the initial boundary value
problem for an abstract integro-differential equation differs substantially from the approach used by L. Pandolfi in~\cite{21} and R. K. Miller and R. L. Wheeler in \cite{33}, \cite{7}.
Moreover, L. Pandolfi studied solvability in function space on a finite time interval $(0,T)$, whereas we study
solvability in the weighted Sobolev spaces $W_{2,\gamma}^{2}({\mathbb R}_+,A^2)$  on the positive semiaxis ${\mathbb R}_+$.

Our proof of the solvability theorem \ref{T:2},  essentially uses the Hilbert structure of the
space $W_{2,\gamma}^{2}({\mathbb R}_+,A_0^2)$, $L_{2,\gamma}({\mathbb R}_+,H)$  and Paley-Wiener theorem, while
in \cite{33}--\cite{21}, considerations are performed in Banach function spaces consisting of smooth functions on a finite time interval $(0, T)$.

\section{Statement of the problem}
Let $H$ be a separable Hilbert space,  and let $A$ be a self-adjoint positive operator with bounded inverse acting
on $H$. 


Consider the following problem for a second-order integro-differential equation on ${\R}_ +   = (0,\infty )$:

\begin{multline}\label{eq:101}
\frac{{d^2 u(t)}}{{dt^2 }} + A^2 u(t) + Bu(t)-\int_0^t {K(t - s)A^2u(s)ds}-\\-\int_0^t {Q(t - s)Bu(s)ds}  = f(t),\quad t\in \R_+,
\end{multline}
\begin{equation}\label{eq:102}
u(+0)= \varphi _0, \quad u^{(1)}(+0)=\varphi_1,
\end{equation}
where $A$ is a positive self-adjoint operator acting in the separable Hilbert space $H$, ${A^*} = A \geqslant {\kappa _0}$ (${\kappa _0}>0$), having the compact inverse operator, $I$ is identity operator in the separable Hilbert space $H$; operator $B$ is symmetric on the $Dom\left( {A^2 } \right)$, nonnegative,  $\left( {Bx,y} \right) = \left( {x,By} \right)$, $\left( {Bx,x} \right) \geqslant 0$ for arbitrary $x,y \in Dom\left( {A^2 } \right)$ and satisfying the inequality  $\left\| {Bx} \right\| \leqslant \kappa \left\| {A^2 x} \right\|$,  $0<\kappa<1$, $x\in Dom\left( {A^2 } \right)$.

We suppose that the kernels $K(t)$ and $Q(t)$ can be represented in the following form
\begin{equation}\label{eq:3g}
K(t) = \sum\limits_{j = 1}^\infty  {{c_j} e^{ - \gamma _j t} }, \quad Q(t) = \sum\limits_{j = 1}^\infty  {{d_j} e^{ - \gamma _j t} },
\end{equation}
where the coefficients $c_j>0$, $d_j\geqslant 0$, $\gamma_{j+1}>\gamma_j> 0$, $j\in \N$, $ \gamma _j  \to +
\infty$ $(j\to +\infty)$ and moreover we suppose that
\begin{equation}\label{eq:91}
\sum\limits_{j = 1}^\infty  {\frac{{c_j }}{{\gamma _j }}}  <  \infty, \quad \sum\limits_{j = 1}^\infty  {\frac{{d_j }}{{\gamma _j }}}  <  \infty.
\end{equation}
The conditions \eqref{eq:91} means that $K(t), Q(t)\in L_1(\mathbb R_+)$. If in addition to \eqref{eq:91}, conditions
\begin{equation}\label{eq:60g}
K(0)=\sum\limits_{j = 1}^\infty  {c_j }  <  + \infty, \quad Q(0)=\sum\limits_{j = 1}^\infty  {d_j }  <  + \infty.
\end{equation}
hold, then the kernels $K(t)$ and $Q(t)$ belong to the space $W^1_1(\mathbb R_+)$.

Let us denote
$$
A_0^2=A^2+B.
$$
Due to known result (theorem \cite{Kato}, p.361) operator $A_0^2$ is self-adjoint positive operator.
We convert the domain $Dom(A_0^{\beta})$ of the operator $A_0^{\beta}$, $\beta>0$ into a Hilbert space $H_{\beta}$ by
endowing  $Dom(A_0^{\beta})$ with the norm $\|\cdot\|_{\beta}=\|A_0^{\beta}\cdot\|$, which is equivalent to the graph norm
of the operator $A_0^{\beta}$.

We note that operator-function $L(\lambda)$ is the symbol of the equation \eqref{eq:101} and it has the following form
\begin{equation}\label{eq:19}
L(\lambda ) = {\lambda ^2}I + {A^2} + B - \hat K(\lambda ){A^2} - \hat Q(\lambda )B,
\end{equation}
where $\hat K(\lambda )$ and $\hat Q(\lambda )$ are the Laplace transforms of the kernels $K(t)$ and $Q(t)$, having the representations 
\begin{equation}\label{eq:LQ}
\hat K(\lambda )=\sum\limits_{j = 1}^\infty
{\frac{{c_j }}{(\lambda  + \gamma _j )}}, \quad \hat Q(\lambda )=\sum\limits_{j = 1}^\infty
{\frac{{d_j }}{(\lambda  + \gamma _j )}},
\end{equation}
$I$ is identity oprerator in the space $H$.

\begin{remark}
Because of the assumptions imposing on operator $A$ and $B$ operator $A_0$ is invertible and the operators $A^2A_0^{-2}$, $BA_0^{-2}$ are bounded and operator $A_0^{-1}$ is compact.
\end{remark}

By $W_{2,\gamma }^n \left( {{\R}_ +  ,A_0^n } \right)$ we denote the Sobolev space of vector functions on the half-axis
${\R}_ + = (0,\infty )$ taking values in $H$ endowed with the norm
$$
\left\| u \right\|_{W_{2,\gamma }^n ({\R}_+, A_0^n )}  \equiv \left(
{\int_0^\infty  {e^{ - 2\gamma t} \left( {\left\| {u^{(n)} (t)}
\right\|_H^2  + \left\| {A_0^n u(t)} \right\|_H^2 } \right)dt} }
\right)^{1/2}, \quad \gamma  \ge 0.
$$

For more information about the spaces $W_{2,\gamma }^n \left( {{\R}_ +  ,A_0^n } \right)$, see monograph \cite{19},
Chapter~1. For $n = 0$ we set $W_{2,\gamma }^0 \left( {\R_ +  ,A_0^0 } \right) \equiv L_{2,\gamma } \left( {\R_ + ,H}
\right)$, and for $\gamma  = 0$, we write $W_{2,0}^n  = W_2^n$.
\begin{definition}
We say that vector function $u$ is a strong solution of the problem \eqref{eq:101}, \eqref{eq:102}, if it belongs to the
space $W_{2,\gamma}^2 ({\mathbb R}_+,A_0^2)$, where $A_0^2=A^2+B$ for some $\gamma \geqslant 0$, satisfies \eqref{eq:101} almost everywhere on
the half-axis ${\mathbb R}_+$, and satisfies the initial conditions \eqref{eq:102}.
\end{definition}

The following theorem give conditions for problem \eqref{eq:101}, \eqref{eq:102} to
be well-solvable.

\begin{theorem}\label{T:2}
If condition \eqref{eq:60g} holds, $f'(t)\in L_{2,\gamma_0}(\mathbb R_+,H)$ for certain $\gamma_0 \geqslant 0$ and $f(0)=0$, $\varphi_0\in H_2$, $\varphi_1\in  H_1$. Then there exists
$\gamma_1 \geqslant \gamma_0$, that for any $\gamma>\gamma_1$ problem \eqref{eq:101}, \eqref{eq:102} is uniquely solvable in the space $W_{2,\gamma }^2 \left( {\mathbb R_ +,A_0^2 } \right)$, and its solution satisfies the inequality
\begin{equation}\label{eq:92}
\left\| u \right\|_{W_{2,\gamma }^2 \left( {\mathbb R_ +  ,A_0^2 } \right)}
\leqslant d\left(\left\| {f'(t)} \right\|_{L_{2,\gamma } \left( {\mathbb R_ + ,H }
\right)}+\left\|A_0^2\varphi_0\right\|_H+\left\|A_0\varphi_1\right\|_H\right),
\end{equation}
where the constant $d$ does not depend on the vector function $f$ and the vectors $\varphi_0$ and $\varphi_1$.
\end{theorem}

Equation~\eqref{eq:101} is related to applications: if $B\equiv 0$ then it is an abstract
form of the Gurtin-Pipkin integro-differential equation modelling the finite-speed heat propagation in
media with memory.
That integrodifferential equation is deduced in \cite{20}.

Equations of the above type are currently investigated by many authors (see \cite{20}, \cite{26}, \cite{21}, \cite{41}--\cite{18}, \cite{418}, \cite{AFG} and references
therein).

We impose the following assumptions.
\begin{enumerate}
  \item The operator $B$ is identically zero,
  \item The operator $A$ has a compact inverse $A^{-1}$
  \item The real-valued function $K(t)$ setisfies the following assumptions
\begin{equation}\label{eq:111}
\sum\limits_{j=1}^{\infty} c_j <+\infty; \quad \sum\limits_{j=1}^{\infty} \frac{c_j}{\gamma_j} <1.
\end{equation}
\end{enumerate}
Note that a detailed structure of the spectrum of the operator-valued function ${L}(\lambda)$ can be described
under the above assumptions.

Let $\{a_j\}_{j=1}^{\infty}$ be eigenvalues of the operator $A$ $(Ae_j=a_j e_j)$, numbered according to the increasing
order: $0< a_1 < a_2 < \dots < a_n < \dots$; $a_n \rightarrow +\infty$ $(n\rightarrow +\infty)$. The corresponding eigenvectors $\{e_j\}_{j=1}^{\infty}$ form an orthonormal basis of the space $H$.

Now we consider the structure of the spectrum of the operator-valued function ${L}_1 (\lambda)$:
$$
{ L}_1 (\lambda)=\lambda^2 I +A^2 - \hat{K}(\lambda)A^2,
$$
where $\hat{K}(\lambda)$ is the Laplace transform of the function $K$.

In the considered case,~\eqref{eq:101} can be decomposed into a countable set of scalar integro-differential
equations
\begin{equation}\label{eq:222}
u_{n}^{(2)}(t)+a_n^2 u_n(t) - \int_{0}^{t} \sum_{k=1}^{\infty} c_k e^{-\gamma_k (t-s)} a_n^2 u_n(s)\,ds=f_n(t),\,t>0
\end{equation}
where $ u_n(t)=(u(t),e_n)$ and $f_n(t)=(f(t),e_n)$, $n=1,2, \dots$. Those equations are projections~\eqref{eq:101} onto the one-dimensional spaces spanned by vectors $\{e_n\}$.

Using the Laplace transform, we naturally arrive at the countable set of meromorphic functions $l_n (\lambda)$:
\begin{equation}\label{eq:333}
l_n(\lambda)=\lambda^2 +a_n^2 -a_n^2\big(\sum_{k=1}^{\infty}
\frac{c_k}{\lambda+\gamma_k}\big), \, n=1,2, \dots;
\end{equation}
which are symbols of the integrodifferential equations given by~\eqref{eq:222}.

We assume that the following condition is satisfied:
\begin{equation}\label{eq:95}
\mathop {\sup }\limits_k \left\{ {\gamma _k (\gamma _{k + 1}  - \gamma _k )} \right\} =  + \infty
\end{equation}

The spectrum of the operator-valued function ${L_1}(\lambda)$ is described as follows.

\begin{theorem}\label{t:3}
The spectrum of the operator-function ${L}_1(\lambda)$ coincides with the closure of the union
of the sets of zeros for the functions
$\{l_n(\lambda)\}_{n=1}^{\infty}$. The zeros of the meromorphic function $l_n(\lambda)$ form
a countable set of real roots $\{\lambda_{n,k}\}_{n=1}^\infty$, satisfying the inequalities
\begin{equation}\label{eq:54}
\begin{array}{cc}
&-\gamma_k< \lambda_{k,n}<x_k<-\gamma_{k-1} \quad [\gamma_0=0],
\\
&\lambda_{k,n}=x_k + \underline{O}\left(\frac{1}{a_n^2}\right)\quad
 k \in {\mathbb N},\, (a_n \rightarrow +\infty).
\end{array}
\end{equation}
where $x_k$ are the real zeros of the function
\begin{equation}\label{eq:55}
g(\lambda) = 1-\sum \limits_{k=1}^{\infty}\frac{c_k}{\lambda +\gamma_k},
\end{equation}
and a pair of complex-conjugate roots
$\{\lambda_n^{\pm}\}_{n=1}^\infty$,
$\lambda_{n}^+=\overline{\lambda_{n}^-}$ such that
\begin{equation}\label{eq:56}
\lambda_n^{\pm}=\pm i \left(a_n+ \underline{O} \left(\frac{1}{a_n}\right)\right)-\frac{1}{2}\sum_{k=1}^{\infty}c_k  + \underline{O} \left(\frac{1}{a_n^2}\right), \quad (a_n
\rightarrow +\infty).
\end{equation}
Thus, the spectrum $\sigma({L}_1)$ of the operator-valued function
${ L}_1(\lambda)$ is representable as
\begin{equation*}
\sigma({L}_1) \equiv
\overline{\big(\cup_{k=1}^{\infty}\cup_{n=1}^{\infty}\{\lambda_{n_k}\}\big)\cup\big(\cup_{n=1}^{\infty}\lambda_{n}^{\pm}\big)},
\end{equation*}
where $\lim\limits_{n\to\infty} \lambda_{k,n}=x_k$, $k=1,2, \dots$.
\end{theorem}

The proof of Theorem~\ref{t:3} is given in \cite{18}. Picture of the spectrum of operator-function $L_1(\lambda)$ is given at Figure 1.

\begin{figure}[h]
      \noindent\centering{
        \small
        \unitlength 1em
         \begin{picture}(25,18)
        \put(0,0){\includegraphics[width=27em,height=18em]{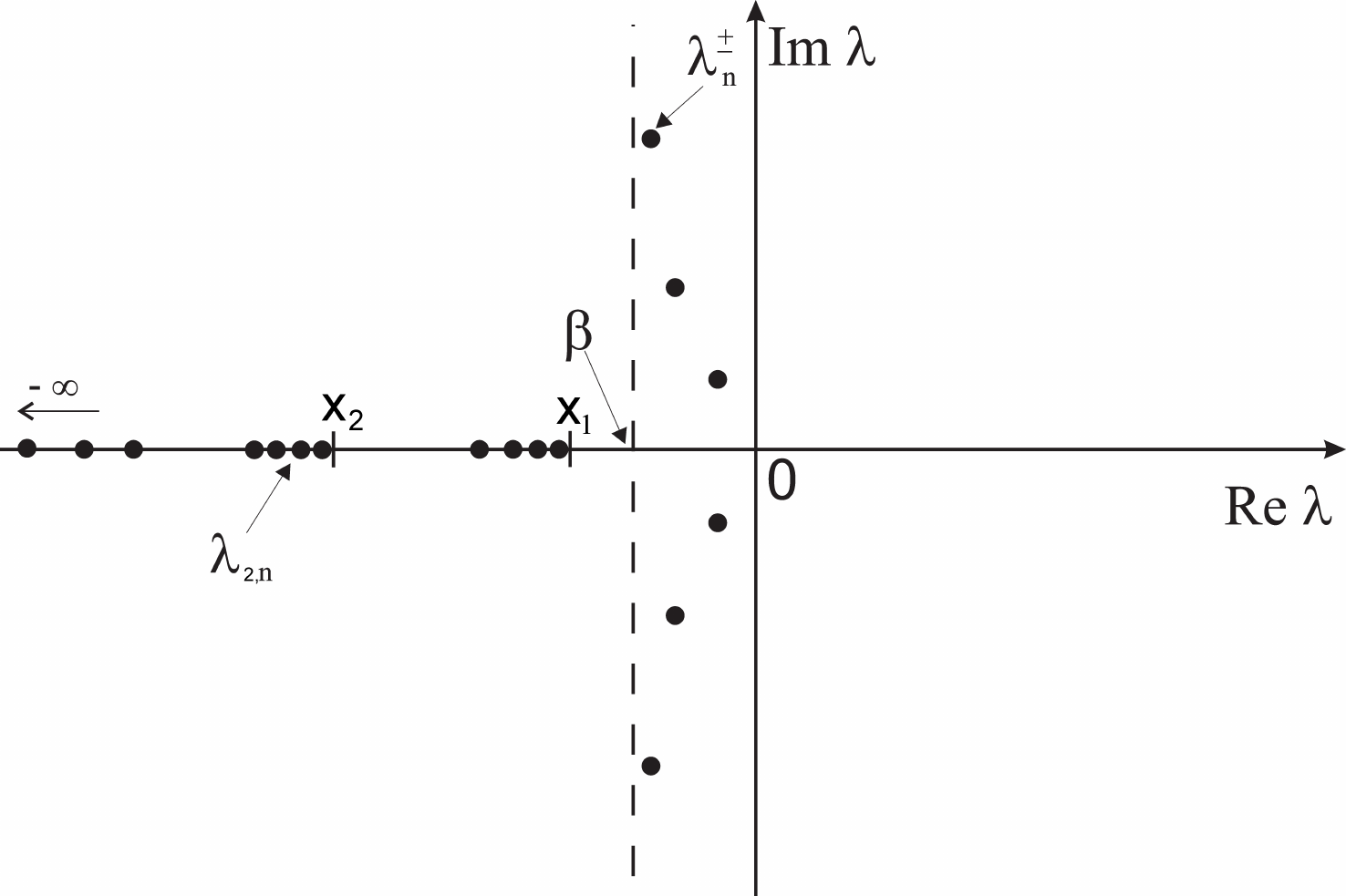}}
       \end{picture}\\}
      \caption{Spectral structure in case  $K(t)\in W^1_1(\R_+)$, $\beta=-\frac{1}{2}K(0)$}
       \label{f:cross_intr}
        \end{figure}

\begin{remark} The spectrum of the operator-valued function ${L_1}(\lambda)$ is located in the left-hand semiplane $\{\lambda: \Re \lambda <0\}$, if

$$
\sum_{j=1}^{\infty}\frac{c_j}{\gamma_j}<1.
$$
If

$$
\sum_{j=1}^{\infty}\frac{c_j}{\gamma_j}>1,
$$
then the accumulation point $x_1$ of the poles, which are eigenvalues $\{\lambda_{1n}\}_{n=1}^{\infty}$ of the operator-valued
function ${L_1}$ is located in the right-hand semi-plane $\{\lambda: \Re \lambda >0\}$; this corresponds to the instability
case.
\end{remark}

We obtain representations of solutions of problem \eqref{eq:101}, \eqref{eq:101} for  $B\equiv 0$ in the form of series in exponentials
corresponding to points of spectrum of the operator-valued function ${L}_1 (\lambda)$.
\begin{theorem}\label{T:3}
Suppose that $f(t)=0$ for $t\in \R_+$, the vector function $u(t)\in W_{2,\gamma }^2 \left( {\R_ + ,A^2 } \right)$, $\gamma>0$ is a strong solution of problem \eqref{eq:101}, \eqref{eq:102} and conditions \eqref{eq:111}, \eqref{eq:95} are satisfied. Then, for each $t\in \R_+$ the solution $u(t)$ of  problem \eqref{eq:101}, \eqref{eq:102} can be represented as the following sum of series:
\begin{multline}\label{eq:uu}
u(t) =\sum\limits_{n = 1}^{\infty}  {\frac{{ \left(\varphi_{1n}+\lambda _n^ +  \varphi_{0n}\right)e^{\lambda _n^ +  t} e_n }}{{l'_n(\lambda _n^ + ) }}}  + \sum\limits_{n = 1}^{\infty}  {\frac{{  \left(\varphi_{1n}+\lambda _n^ -  \varphi_{0n}\right) e^{\lambda _n^ -  t} e_n }}{{l'_n(\lambda _n^ - )}}}  + \\+\sum\limits_{n = 1}^{\infty}  {\left( {\sum\limits_{k = 1}^{\infty} {\frac{{ \left(\varphi_{1n}+\lambda _{k,n}  \varphi_{0n}\right)e^{\lambda _{k,n} t} }}{{l'_n(\lambda _{k,n}) }}} } \right)} e_n, 
\end{multline}
convergent in the norm of $H$, where $\varphi_{in}=(\varphi_{i},e_n)$, $(i=0,1)$, 
$\lambda _{k,n}$ are the real zeros of the meromorphic function  $l_n (\lambda )$, which satisfy the inequalities \eqref{eq:54}
and $\lambda _n^ \pm$ is a pair of complex conjugate zeros, $\lambda _n^ + =
\overline {\lambda _n^ -  }$, asymptotically representable in the form \eqref{eq:56}.
\end{theorem}

\begin{theorem}\label{T:4}
If $\varphi_0\in H_2$ and $\varphi_1\in H_1$, 
then the series obtained from \eqref{eq:uu} by $p$-fold termwise differentiation with respect to $t$ for $p=0,1,2$ converges in the space $H_{2-p}$ uniformly with respect
to $t$ on any interval $[t_0,T]$, where $0<t_0<T<+\infty$. Moreover, for all $t\in [t_0,T]$ one has the
estimates
$$
\left\| {\sum\limits_{n = 1}^\infty  {u_n^{(p)} (t)e_n } } \right\|_{H_{2 - p} }^2 \leqslant d(\|A\varphi _{1}\|^2+\|A^2\varphi _{0}\|^2), \quad p=0,1,2,
$$
where $u_n(t)=(u(t),e_n)$, with constant $d$ independent of the vector functions $\varphi _{1}$ and $\varphi _{0}$.

Further, if there are finitely many summands in \eqref{eq:222} (i.e., $c_j=0$ for $j>N$, $N\in \mathbb N$), then one can
set $t_0=0$.
\end{theorem}
\begin{theorem}\label{T:5}
Suppose that $f(t)\in C\left( {[0,T] , H} \right)$ for arbitrary $T>0$ and suppose that the vector-function $u(t)\in W_{2,\gamma }^2 \left( {\R_ + ,A^2 } \right)$ for some $\gamma>0$, is a strong solution of problem \eqref{eq:101}, \eqref{eq:102}, conditions \eqref{eq:111}, \eqref{eq:95} are satisfied, and $\varphi_0=\varphi_1=0$.  Then, for each $t\in \R_+$, the solution $u(t)$ of problem \eqref{eq:101}, \eqref{eq:102} can be
represented as the following sum of series: 
\begin{equation}\label{eq:u1}
u(t) = \sum\limits_{n = 1}^\infty  {\omega_ n (t,\lambda _n^ +  )e_n }  + \sum\limits_{n = 1}^\infty  {\omega_n  (t,\lambda _n^ -  )e_n }  + \sum\limits_{n = 1}^\infty  {\left( {\sum\limits_{k = 1}^\infty  {\omega_n  (t,\lambda _{k,n}  )} } \right)e_n },
\end{equation}
converging in the norm of $H$, where 
$$
\omega_n (t,\lambda ) = \frac{{ {\int\limits_0^t {f_n(\tau )e^{\lambda (t-\tau)}d\tau } } }}
{{l'_n (\lambda )}},
$$
$f_n(t)=(f(t),e_n)$.
\end{theorem}
\begin{theorem}\label{T:6}
Let the assumptions of Theorem  \ref{T:5} be satisfied together with the condition $\sum\limits_{j = 1}^\infty  {{{\gamma _j^{-3/2} }}}  < \infty$. Then the series obtained from \eqref{eq:u1} by $p$-fold term-by-term differentiation with respect to $t$ for $p=0,1,2$ is convergent in the space $H_{2-p}$ uniformly with respect to $t$ on any interval $[t_0,T]$, where $t_0<T<+\infty$ ($t_0=0$ for $p=0,1$ and $t_0>0$ for $p=2$); moreover
the following estimates hold for all $t\in [t_0,T]$:
\begin{equation}\label{eq:e1}
\left\| {\sum\limits_{n = 1}^\infty  {u_n (t)e_n } } \right\|_{H_2 }^2 \leqslant d\left\| {A^2 f(t)} \right\|_{L_{2,\gamma } (\R_ +  ,H)}, 
\end{equation}
\begin{multline}\label{eq:e2}
\left\| {\sum\limits_{n = 1}^\infty  {u_n^{(p)} (t)e_n } } \right\|_{H_{2 - p} }^2  \leqslant \\\leqslant d\left( {\left\| {A^{2 - p} f^{(p)} (t)} \right\|_{L_{2,\gamma } (\R_ +  ,H)}^2  + \left\| {A^{2 - p}f(0)} \right\|_H^2  + (p - 1)\left\| {f'(0)} \right\|_H^2 } \right), \\\quad p=1,2,
\end{multline}
where $u_n(t)=(u(t),e_n)$.
\end{theorem}

The series \eqref{eq:uu} and \eqref{eq:u1} are obtained by applying the inverse Laplace transform to the solution
of problem \eqref{eq:101}, \eqref{eq:102} with the use of integration over rectangular contours separating the points $\gamma_j$. (The construction of these contours is described in \cite{412}.) Here the estimates of the operator function $L_1^{-1}(\lambda)$ on these contours play an important role. Besides, the proof of Theorems \ref{T:3} -- \ref{T:6} heavily uses the representations \eqref{eq:54} and \eqref{eq:56} obtained in \cite{412}. The proofs of Theorems \ref{T:3} -- \ref{T:6} are obtained in \cite{OT}.

Throughout the paper, the expression $D\lesssim E$ stands for the inequality $D\leqslant const\cdot E,$ with a positive
constant $const$; the expression $D\approx E$ means that $D \lesssim E\lesssim 
D.$  The symbols := and =: are used for the
introduction of new quantities.

\section{Proofs.}
We begin the proof of the theorem \ref{T:2} in the case of homogeneous (zero) initial conditions ($\varphi_0=\varphi_1=0$).  We use Laplace transformation in oder to prove the correct solvability of the problem \eqref{eq:101} and \eqref{eq:102}. Now we are going to remind the base assertions that will be used  later.

\begin{definition}
We denote by $H_2 (\Re\lambda > \gamma ,H)$ the Hardy space of vector-functions $\hat f(\lambda)$ taking
values in the space $H$, holomorphic (analytic) in the semiplane $\{\lambda\in\mathbb C:\Re
\lambda>\gamma \geqslant 0\}$ endowed with the norm
\begin{equation}\label{eq:15}
\mathop {\sup }\limits_{x > \gamma } \int_{ - \infty }^{ + \infty }
{\left\| {\hat f(x + iy)} \right\|} _H^2 dy < \infty, \quad
(\lambda=x+iy).
\end{equation}
\end{definition}
We formulate well-known Paley-Wiener theorem for Hardy space   $H_2 (\Re \lambda > \gamma ,H)$.

\bigskip
{\noindent\bf Theorem} (Paley-Wiener). {\it 1) The space  $H_2 (\Re
\lambda > \gamma ,H)$ coincides with the set of vector-functions (Laplace transformations) representing in the form 
\begin{equation}\label{eq:16}
\hat f(\lambda ) = \frac{1}{{\sqrt {2\pi } }}\int_0^\infty  {e^{ -
\lambda t} f(t)dt},
\end{equation}
with   vector-function $f(t)\in L_{2,\gamma}\left(\mathbb R_+,H\right)$, $\lambda \in \mathbb C$,
$\Re \lambda >\gamma \geqslant 0$.

\noindent 2) There exists unique  vector-function  $f(t)\in L_{2,\gamma}\left(\mathbb R_+,H\right)$ for arbitrary vector-function $\hat f(\lambda )\in H_2 (\Re\lambda> \gamma ,H)$ and the following inversion formula take place 
\begin{equation}\label{eq:17}
f(t) = \frac{1}{{\sqrt {2\pi } }}\int_{ - \infty }^{ + \infty }
{\hat f(\gamma  + iy)e^{(\gamma  + iy)t} dy}, \quad t\in \R_+,\quad
\gamma \geqslant 0
\end{equation}

\noindent 3)  The following equality take place for vector-function  $\hat f(\lambda )\in H_2 (\Re
\lambda > \gamma ,H)$ and  $f(t)\in L_{2,\gamma}\left(\R_+,H\right)$,
connected by relation \eqref{eq:16}:}
\begin{multline}\label{eq:18}
\|\hat f\|^2_{H_2(\Re \lambda>\gamma,H)}\equiv \mathop {\sup
}\limits_{x> \gamma } \int_{ - \infty }^{ + \infty } {\left\| {\hat f(x + iy)}
\right\|_H^2 } dy =\\= \int_{0}^{ + \infty } {e^{ - 2\gamma t}
\left\| {f(t)} \right\|_H^2 dt}\equiv \|f\|_{L_{2,\gamma}\left(\R_+,H\right)}^2
\end{multline}

\bigskip
\noindent The theorem formulated above is well-known for the scalar functions. However it is easily generalized for the vector-functions taking values in the separable Hilbert space.(see ~\cite{32}).

{\bf Proof of the theorem \ref{T:2}.}  Let us consider the case $f'(t)\in L_{2,\gamma_2 } \left( {\R_ + ,H }\right)$, $f(0)=0$. The Laplace transformation of the strong solution of the problem \eqref{eq:101} and \eqref{eq:102} with initial conditions $u(+0)=0$, $u^{(1)}(+0)=0$ has the following representation 
\begin{equation}\label{eq:501}
\hat u(\lambda ) = L^{ - 1} (\lambda )\hat f(\lambda ).
\end{equation}
Here the operator-function $L(\lambda)$ is the symbol of the equation \eqref{eq:101} and it has the following form
\begin{equation}\label{eq:19*}
L(\lambda ) = {\lambda ^2}I + {A^2} + B - \hat K(\lambda ){A^2} - \hat Q(\lambda )B,
\end{equation}
where $\hat K(\lambda )$ and $\hat Q(\lambda )$ are the Laplace transforms of the kernels $K(t)$ and $Q(t)$, having the representations 
\begin{equation}\label{eq:LQ*}
\hat K(\lambda )=\sum\limits_{j = 1}^\infty
{\frac{{c_j }}{(\lambda  + \gamma _j )}}, \quad \hat Q(\lambda )=\sum\limits_{j = 1}^\infty
{\frac{{d_j }}{(\lambda  + \gamma _j )}},
\end{equation}
$\hat f(\lambda)$ is Laplace transform of vector-valued function $f(t)$, $I$ is identity oprerator in the space $H$.

We suppose there exists $\gamma ^*  \ge 0$, so that $u \in W_{2,\gamma ^* }^2 \left( {\R_ + ,A_0^2} \right)$. We need this supposition in oder to use the Laplace transformation to the equation \eqref{eq:101}. 

It is sufficient (for proof of the theorem \ref{T:2}) to prove that vector-functions   
 $\lambda^2 \hat u(\lambda )$ and $A_0^2\hat u(\lambda )$ belong to Hardy space
 $H_2 (\Re\lambda > \gamma ,H)$ for some $\gamma>\gamma_0\ge 0$. Then we shall obtain by the Paley-Wiener theorem that vector-functions $\frac{{\displaystyle d^2 u}}{{\displaystyle dt^2 }}$ and $A_0^2u(t)$ belong to the space $L_{2,\gamma } \left( {\R_ +  ,H} \right)$ and hence $u(t) \in W_{2,\gamma}^2 \left( {\R_ + ,A_0^2 }
\right)$. Thus we shall prove the solvability of the problem  \eqref{eq:101} and \eqref{eq:102} with homogeneous initial conditions in the space  $W_{2,\gamma}^2 \left( {\R_ + ,A_0^2 } \right)$.

Further we shall use the following Lemma.
\begin{lemma}\label{L1}
There exists such $\gamma>0$ that the inequality
\begin{equation}\label{eq:Lemma1}
\mathop {\sup }\limits_{\operatorname{Re} \lambda  > \gamma  > 0} \left\| {\frac{1}{\lambda }A_0^2{{\left( {{\lambda ^2}I + A_0^2} \right)}^{ - 1}}} \right\| < \frac{{const}}{\tau }, \quad \lambda=\tau+i\nu
\end{equation}
is valid.
\end{lemma}

{\bf Proof of the Lemma \ref{L1}.} Taking into account that operator $A_0$ is self-adjoint we shall use the spectral theorem (see \cite{Kato}, pp. 452-453). Suppose that $\lambda =\tau+i\nu$, $(\tau, \nu \in \mathbb R)$ and 
$a\in \sigma(A_0)$ that is $a$ belongs to the spectra of operator $A_0$. Due to the spectral theorem it is sufficient to prove the following estimate 
\begin{equation}\label{est1}
\frac{{{a^2}}}{{{{({\tau ^2} + {\nu ^2})}^{1/2}}{{(({\tau ^2} - {\nu ^2}) + {a^2})}^2} + 4{\tau ^2}{\nu ^2}{)^{1/2}}}} \leqslant  \frac{{const}}{\tau }, \quad \tau  \geqslant \gamma  > 0.
\end{equation}

In order to do this let us estimate function 
$$
f(a,\tau ,\nu ) = ({\tau ^2} + {\nu ^2}){(({\tau ^2} - {\nu ^2}) + {a^2})^2} + 4{\tau ^2}{\nu ^2})
$$
from the below. Let us consider the constant $d \in (0,1)$. Then the following estimate is valid
\begin{multline*}
f(a,\tau ,\nu ) \geqslant \min \left\{ {\mathop {\min }\limits_{{\nu ^2} \in [0,d{a^2}]} f(a,\tau ,\nu ),\mathop {\min }\limits_{{\nu ^2} \in [d{a^2}, + \infty ]} f(a,\tau ,\nu )} \right\} \geqslant \\ \geqslant \min \left\{ {({\tau ^2} + (1 - d){a^2})^2{\tau ^2},({\tau ^2} + d{a^2})4d{a^2}{\tau ^2}} \right\}.
\end{multline*}
Hence we have 
$$
\frac{{{a^2}}}{{{{(f(a,\tau ,\nu ))}^{1/2}}}} \leqslant {a^2}\max \left[ {\frac{1}{{{{(\tau^2 + (1 - d){a^2})}^{1/2}}\tau }},\frac{1}{{{{({\tau ^2} + d{a^2})}^{1/2}}2\sqrt d a\tau }}} \right] \leqslant 
$$
\begin{multline}\label{est2}
\leqslant \max \left[ {\frac{1}{{\tau {{\left( {\frac{{{\tau ^2}}}{{{a^2}}} + (1 - d)} \right)}^{1/2}}}},\frac{1}{{2\sqrt d \tau {{\left( {\frac{{{\tau ^2}}}{{{a^2}}} + d)} \right)}^{1/2}}}}} \right] \leqslant \max \left[ {\frac{1}{{\tau (1 - d)}},\frac{1}{{2\tau d}}} \right].
\end{multline}
Puting $d=1/3$ we obtain from \eqref{est2} the desired inequality \eqref{est1}. Lemma \ref{L1} is proved.

\begin{lemma}\label{L2}
Suppose the conditions of the theorem 1 are satisfies. Then there exists $\gamma>0$ such that operator-valued function $\left(I-V(\lambda)\right)^{-1}$, where 
\begin{equation}\label{V}
V(\lambda ) = \hat K(\lambda ){A^2}{\left( {{\lambda ^2}I + A_0^2} \right)^{ - 1}} + \hat Q(\lambda )B{\left( {{\lambda ^2}I + A_0^2} \right)^{ - 1}}
\end{equation}
is analytic at half-plane $\left\{ {\lambda :\operatorname{Re} \lambda  > \gamma } \right\}$ and satisfies the following inequality
\begin{equation}\label{est3}
\mathop {\sup }\limits_{\lambda :\operatorname{Re} \lambda  > \gamma } \left\| {{{\left( {I - V(\lambda )} \right)}^{ - 1}}} \right\| \leqslant const.
\end{equation}
\end{lemma}

{\bf Proof of the Lemma \ref{L2}.} On the base of the conditions \eqref{eq:60g} it is easy to check that 
\begin{equation}\label{est4}
\mathop {\sup }\limits_{\operatorname{Re} \lambda  > 0} \left| {\hat K(\lambda )} \right| = \mathop {\sup }\limits_{\operatorname{Re} \lambda  > 0} \left| {\sum\limits_{j = 1}^\infty  {\frac{{{c_j}}}{{(\lambda  + {\gamma _j})}}} } \right| \leqslant \left( \sum\limits_{j = 1}^\infty  {c_j} \right)\mathop {\sup }\limits_{\operatorname{Re} \lambda  > 0} \frac{1}{{\left| {\lambda  + {\gamma _1}} \right|}} \leqslant \frac{{{K(0)}}}{{|\lambda |}},
\end{equation}
\begin{equation}\label{est5}
\mathop {\sup }\limits_{\operatorname{Re} \lambda  > 0} \left| {\hat Q(\lambda )} \right| = \mathop {\sup }\limits_{\operatorname{Re} \lambda  > 0} \left| {\sum\limits_{j = 1}^\infty  {\frac{{{d_j}}}{{(\lambda  + {\gamma _j})}}} } \right| \leqslant \left( \sum\limits_{j = 1}^\infty  {d_j}  \right)\mathop {\sup }\limits_{\operatorname{Re} \lambda  > 0} \frac{1}{{\left| {\lambda  + {\gamma _1}} \right|}} \leqslant \frac{{{Q(0)}}}{{|\lambda |}}.
\end{equation}
Then using Lemma 1 we obtain that there exists such $\gamma_1>0$ that for $\operatorname{Re} \lambda  > \gamma_1$ we have
\begin{multline}\label{est6}
\left\| {\hat K(\lambda ){A^2}{{\left( {{\lambda ^2}I + A_0^2} \right)}^{ - 1}}} \right\| \leqslant \left| {\hat K(\lambda )} \right|\left\| {{A^2}A_0^{ - 2}} \right\|\left\| {A_0^2{{\left( {{\lambda ^2}I + A_0^2} \right)}^{ - 1}}} \right\| \leqslant \\ \leqslant \left|K(0)\right|\left\| {{A^2}A_0^{ - 2}} \right\|\left\| {\frac{1}{\lambda }A_0^2{{\left( {{\lambda ^2}I + A_0^2} \right)}^{ - 1}}} \right\| \leqslant \frac{{const}}{{\operatorname{Re} \lambda }}.
\end{multline}
\begin{multline}\label{est6*}
\left\| {\hat Q(\lambda ){B}{{\left( {{\lambda ^2}I + A_0^2} \right)}^{ - 1}}} \right\| \leqslant \left| {\hat Q(\lambda )} \right|\left\| {{B}A_0^{ - 2}} \right\|\left\| {A_0^2{{\left( {{\lambda ^2}I + A_0^2} \right)}^{ - 1}}} \right\| \leqslant \\ \leqslant \left|Q(0)\right|\left\| {{B}A_0^{ - 2}} \right\|\left\| {\frac{1}{\lambda }A_0^2{{\left( {{\lambda ^2}I + A_0^2} \right)}^{ - 1}}} \right\| \leqslant \frac{{const}}{{\operatorname{Re} \lambda }}.
\end{multline}

Uniting the estimates  \eqref{est6}, \eqref{est6*} we have
\begin{equation}\label{est7}
\left\| {V(\lambda )} \right\| \leqslant \frac{{const}}{{\operatorname{Re} \lambda }}.
\end{equation}
and thus it is possible to choose such $\gamma>0$ that 
$$
\mathop {\sup }\limits_{\operatorname{Re} \lambda  > {\gamma}} \left\| {V(\lambda )} \right\| < 1.
$$
Thus we proved Lemma \ref{L2}.

\begin{lemma}\label{L3}
The following estimate 
\begin{equation}\label{est8}
\left\| {\lambda {{\left( {{\lambda ^2}I + A_0^2} \right)}^{ - 1}}} \right\| \leqslant \frac{1}{{\left| {\operatorname{Re} \lambda } \right|}}
\end{equation}
\end{lemma}
is valid.

{\bf Proof of the Lemma \ref{L3}.} Suppose parameter $a$ belongs to the spectra of operator $A$. Using the spectral theorem it is sufficient to note that the following chain of inequalities
$$
\left| {\lambda {{\left( {{\lambda ^2} + {a^2}} \right)}^{ - 1}}} \right| = \frac{{{{({\tau ^2} + {\nu ^2})}^{1/2}}}}{{{{({\tau ^2} + {{(\nu  + a)}^2})}^{1/2}}{{({\tau ^2} + {{(\nu  - a)}^2})}^{1/2}}}} \leqslant \frac{1}{{|\tau |}}, \quad \lambda=\tau+i\nu.
$$
is valid. Thus the estimate \eqref{est8}. Lemma \ref{L3} is proved.

Now using Lemmas \ref{L1}-\ref{L3} and the representation 
\begin{multline}\label{rep1}
\hat u(\lambda ) = \frac{1}{\lambda }\left[{{\left( {{\lambda ^2}I + A_0^2} \right)}^{ - 1}}\left( {I - \hat K(\lambda ){A^2}{{\left( {{\lambda ^2}I + A_0^2} \right)}^{ - 1}} -} \right. \right.\\\left.\left.{-\hat Q(\lambda )B{{\left( {{\lambda ^2}I + A_0^2} \right)}^{ - 1}}} \right)^{ - 1}\lambda f(\lambda )\right]
\end{multline}
We shall prove that vector-valued functions ${\lambda ^2}\hat u(\lambda )$ and $A_0^2\hat u(\lambda )$ belong to Hardy space ${H_2}\left( {\operatorname{Re} \lambda  > \gamma ,H} \right)$.
Owing to the fact that $f'(t) \in {L_{2,\gamma }}\left( {{\mathbb{R}_ + },H} \right)$ and $f(0)=0$ we obtain that vector function  ${\lambda}\hat f(\lambda )$ belongs to Hardy space ${H_2}\left( {\operatorname{Re} \lambda  > \gamma_0 ,H} \right)$. From Lemmas 1, 2 have we that operator-valued function 
\begin{multline*}
A_0^2\frac{1}{\lambda }\left[{{\left( {{\lambda ^2}I + A_0^2} \right)}^{ - 1}}\left( {I - \hat K(\lambda ){A^2}{{\left( {{\lambda ^2}I + A_0^2} \right)}^{ - 1}} -\hat Q(\lambda )B{{\left( {{\lambda ^2}I + A_0^2} \right)}^{ - 1}}} \right)^{ - 1}\right]
\end{multline*}
is bounded and analytic in the half-plane $\left\{ {\lambda :\operatorname{Re} \lambda  > \gamma } \right\}$. Let us choose ${\gamma _1} = \max ({\gamma _0},\gamma )$. Hence vector-valued function $A_0^2\hat u(\lambda )\in {H_2}\left( {\operatorname{Re} \lambda  > \gamma_1 ,H} \right)$. Moreover we obtain the following estimate
\begin{multline}\label{est9}
{\left\| {A_0^2\hat u(\lambda )} \right\|_{{H_2}\left( {\operatorname{Re} \lambda  > \gamma_1 ,H} \right)}} \leqslant \mathop {\sup }\limits_{\operatorname{Re} \lambda  > \gamma_1 } \left\| {A_0^2\frac{1}{\lambda }{{\left( {{\lambda ^2}I + A_0^2} \right)}^{ - 1}}{{\left( {I - V(\lambda )} \right)}^{ - 1}}} \right\| \cdot \\ \cdot {\left\| {\lambda f(\lambda )} \right\|_{{H_2}\left( {\operatorname{Re} \lambda  > \gamma_1 ,H} \right)}}.
\end{multline}
Using the Paley-Wiener theorem we have 
\begin{equation}\label{est10}
{\left\| {A_0^2u(t)} \right\|_{{L_{2,\gamma_1 }}\left( {{\mathbb{R}_ + },H} \right)}} \leqslant const{\left\| {f'(t)} \right\|_{{L_{2,\gamma_1 }}\left( {{\mathbb{R}_ + },H} \right)}}.
\end{equation}

In turn from Lemmas \ref{L2} and \ref{L3} we obtain that operator-function $\lambda {\left( {{\lambda ^2}I + A_0^2} \right)^{ - 1}}{\left( {I - V(\lambda )} \right)^{ - 1}}$ is bounded and analytic at half-plane $\left\{ {\lambda :\operatorname{Re} \lambda  > \gamma } \right\}$. Thus vector-function ${\lambda ^2}\hat u(\lambda ) \in {H_2}\left( {\operatorname{Re} \lambda  > \gamma_1 ,H} \right)$ and moreover the following estimate is valid
\begin{multline}\label{est11}
{\left\| {\lambda^2\hat u(\lambda )} \right\|_{{H_2}\left( {\operatorname{Re} \lambda  > \gamma_1 ,H} \right)}} \leqslant \mathop {\sup }\limits_{\operatorname{Re} \lambda  > \gamma_1 } \left\| {\lambda {{\left( {{\lambda ^2}I + A_0^2} \right)}^{ - 1}}{{\left( {I - V(\lambda )} \right)}^{ - 1}}} \right\| \cdot \\ \cdot {\left\| {\lambda f(\lambda )} \right\|_{{H_2}\left( {\operatorname{Re} \lambda  > \gamma_1 ,H} \right)}}.
\end{multline}

Using the Paley-Wiener theorem we have 
\begin{equation}\label{est12}
{\left\| {\frac{{{d^2}u}}{{d{t^2}}}} \right\|_{{L_{2,\gamma_1 }}\left( {{\mathbb{R}_ + },H} \right)}} \leqslant const{\left\| {f'(t)} \right\|_{{L_{2,\gamma_1 }}\left( {{\mathbb{R}_ + },H} \right)}}.
\end{equation}

On the base of inequalities \eqref{est10} and \eqref{est12} we obtain the estimate
\begin{equation}\label{est13}
{\left\| {u(t)} \right\|_{W_{_{2,\gamma_1 }}^2\left( {{\mathbb{R}_ + },A_0^2} \right)}} \leqslant const{\left\| {f'(t)} \right\|_{{L_{2,\gamma_1 }}\left( {{\mathbb{R}_ + },H} \right)}}.
\end{equation}

Thus we proved theorem \ref{T:2} for homogeneous initial data \eqref{eq:102}.

Let us consider now the general case that is the problem \eqref{eq:101}, \eqref{eq:102} with nonzero initial data $\varphi_0$ and $\varphi_1$. We will look for the solution of the problem  \eqref{eq:101}, \eqref{eq:102} in the form
$$
u(t) = \cos (A_0t)\varphi _0  + (A_0)^{ - 1} \sin (A_0t)\varphi _1 +\omega(t),
$$
with new unknown vector-function $\omega (t)$. Then for vector-function $\omega (t)$ we obtain the problem with homogeneous initial conditions
\begin{equation}\label{eq:101*}
\frac{{d^2 \omega(t)}}{{dt^2 }} + A_0^2 \omega(t) -\int_0^t {K(t -s)A^2 \omega(s)ds}-\int_0^t {Q(t -s)B \omega(s)ds} = f_1(t),\quad t\in \R_+,
\end{equation}
\begin{equation}\label{eq:102*}
\omega(+0)=\omega^{(1)}(+0)=0.
\end{equation}
where $f_1(t)=f(t)+h(t)$ and vector-function 
\begin{multline}\label{eq:h}
h(t) = \int_0^t {K(t - s)A^2\left( {\cos (A_0s)\varphi _0  + (A_0)^{ - 1}
\sin (A_0s)\varphi _1 } \right)ds}+\\+\int_0^t {Q(t - s)B\left( {\cos (A_0s)\varphi _0  + (A_0)^{ - 1}
\sin (A_0s)\varphi _1 } \right)ds}.
\end{multline}
It is sufficient for us to prove that $h'(t) \in L_{2,\gamma_0 } \left( {\mathbb R_+ ,H} \right)$ and $h(0)=0$. Using integration by parts we obtain the following representation for vector-function 
\begin{equation}\label{eq:h'}
h'(t)=g_1(t)+A^2g_2(t)+Bg_3(t),
\end{equation}
where
$$
g_1(t) = K(t)A^2\varphi _0+Q(t)B\varphi _0,
$$
\begin{equation}\label{eq:g2}
g_2(t)=\int_0^t {K(t - s)\left( -A_0{\sin (A_0s)\varphi _0  + \cos (A_0s)\varphi _1 } \right)ds},
\end{equation}
$$
g_3(t)=\int_0^t {Q(t - s)\left( -A_0{\sin (A_0s)\varphi _0  + \cos (A_0s)\varphi _1 } \right)ds}.
$$
Integrating by parts we have
\begin{multline*}
\int_0^t {e^{ - \gamma _j (t - s)} \cos (A_0s)ds}  =\\=\left( {A_0^2  + \gamma _j^2I } \right)^{ - 1} \left( {\gamma _j \left( {\cos (A_0t) -  e^{ - \gamma _j t} I} \right) + A_0\sin (A_0t)} \right),
\end{multline*}
\begin{multline}\label{eq:int}
\int_0^t {e^{ - \gamma _j (t - s)} \sin (A_0s)ds}  = \\=\left( {A_0^2  + \gamma _j^2 I} \right)^{ - 1} \left( {A_0\left( {e^{ - \gamma _j t}I  - \cos (A_0t)} \right) + \gamma _j \sin (A_0t)} \right).
\end{multline}
We need the following proposition 
\begin{proposition}\label{n:4}
The inequality
\begin{equation}\label{eq:prop}
\left\| {\left( {A_0^2  + \gamma _j^2 I} \right)^{ - 1} } \right\|_H^2  \lesssim \gamma _j^{ - 2} \left\| {A_0^{ - 1} } \right\|_H^2, \quad j\in \mathbb N
\end{equation}
is valid.
\end{proposition}
Really for arbitrary vector $x\in H$ such that $\|x\|_H=1$ the following chain of inequalities is true 
\begin{multline*}
\left\| {\left( {A_0^2  + \gamma _j^2 I} \right)^{ - 1} x} \right\|_H^2  = \sum\limits_{n = 1}^\infty  {\left( {a_n^2  + \gamma _j^2 } \right)^{ - 2} \left| {x_n } \right|^2 }  \lesssim \\\lesssim \sum\limits_{n = 1}^\infty  {\left( {\gamma _j a_n } \right)^{ - 2} \left| {x_n } \right|^2 }  = \left\| {\gamma _j^{ - 1} A_0^{ - 1} x} \right\|_H^{2},
\end{multline*}
where $x_n=(x,e_n)$, $n\in \N$ and $ \{ e_j \} _{j = 1}^\infty$ the orthonormal basis composed of the eigenvectors of the  self-adjoint operator operator $A_0$ corresponding to eigenvalues $a_j$, i.e.,such that $ A_0e_j = a_je_j$  for  $j\in {\N}$. The eigenvalues $a_j$ are numbered in increasing order with multiplicity taken into account: $0 < a_1  \leqslant a_2 \leqslant ...$; $a_n\to +\infty$ as $n\to +\infty$. 

In what follows we will use the following well-known estimates: 
\begin{equation}\label{eq:cos}
\|\cos (A_0t)\|_H\leqslant1, \quad \|\sin (A_0t)\|_H\leqslant1, \quad t\in \mathbb R_+.
\end{equation}
Let us suppose that condition  \eqref{eq:60g} is satisfied and initial vectors $\varphi_0\in D(A_0^2)$, $\varphi_1\in D(A_0)$.
On the base of \eqref{eq:int} let us estimate the norms $\left\| {g_1(t)} \right\|_{L_{2,\gamma } \left( {\R_+ ,H} \right)}$,  $\left\| A^2g_2(t)\right\|_{L_{2,\gamma } \left( {\R_+ ,H}\right)}$, $\left\|Bg_3(t) \right\|_{L_{2,\gamma } \left( {\R_+ ,H} \right)}$. We will use the following remark.
\begin{remark}\label{r1}
From the definition of the operator $A_0^2$ the evident inequalities 
$$
\left\| A^2x\right\|_{L_{2,\gamma } \left( {\R_+ ,H}\right)}\leqslant \left\| A_0^2x\right\|_{L_{2,\gamma } \left( {\R_+ ,H}\right)}, \quad x\in Dom(A_0^2) 
$$
\begin{equation}\label{eq:**}
\left\|Bx\right\|_{L_{2,\gamma } \left( {\R_+ ,H}\right)}\leqslant \left\| A_0^2x\right\|_{L_{2,\gamma } \left( {\R_+ ,H}\right)}.
\end{equation}
follows.
\end{remark}
Due to the Remark \ref{r1} estimate the vector-valued functions $A_0^2g_2(t)$ and $A_0^2g_3(t)$. 

Thus we have the following estimate
\begin{multline}\label{eq:g1}
\left\| {g_1(t)} \right\|_{L_{2,\gamma } \left( {\R_+ ,H} \right)}=
\left\| {\sum\limits_{j = 1}^\infty  {c_j e^{ - \gamma _j t} } A^2 \varphi _0 } +\sum\limits_{j = 1}^\infty  {d_j e^{ - \gamma _j t} } B \varphi _0 \right\|_{L_{2,\gamma } (\R_ +  ,H)}\lesssim \left\| A_0^2\varphi _0 \right\|_H.
\end{multline}
Owing to \eqref{eq:int} vector-valued function $A_0^2g_2(t)$ has the representation
$$
A_0^2g_2(t)=  {\int\limits_0^t {\sum\limits_{j = 1}^\infty  {c_j e^{ - \gamma _j (t - s)} A_0^2\left( { - A_0 \sin \left( { As} \right)\varphi _0  + \cos \left( {A_0s} \right)\varphi _1 } \right)} ds} } =
$$
\begin{multline*}
=\sum\limits_{j = 1}^\infty  c_j A_0^2\left( {A_0^2  + \gamma _j^2 I} \right)^{ - 1} \left( { - A_0} \right)\left(A_0\left( {e^{ - \gamma _j t}I  - \cos \left( {A_0t} \right)} \right) + \gamma _j \sin \left( {A_0t} \right) \right) \varphi _0   + \\+
\sum\limits_{j = 1}^\infty  c_j A_0^2\left( {A_0^2  + \gamma _j^2 I} \right)^{ - 1}  \left( \gamma _j \left( {\cos \left( {A_0t} \right) - e^{ - \gamma _j t} I} \right) +A_0\sin \left( {A_0t} \right) \right) \varphi _1.
\end{multline*}
Further we need the following evident proposition.
\begin{proposition}
The following inequalities 
$$
\left\|A_0^2\left( {A_0^2  + \gamma _j^2 I} \right)^{ - 1}\right\|_H \leqslant const, \quad j\in \mathbb N
$$
are valid.
\end{proposition}
On the base of propositions 3.1, 3.2 and inequalities \eqref{eq:cos} we have the estimates
\begin{multline*}
\left\|A_0^2g_2(t)\right\|_H\leqslant \left\| \sum\limits_{j = 1}^\infty  c_jA_0^2\left( {A_0^2  + \gamma _j^2 I} \right)^{ - 1}A_0^2e^{ - \gamma _j t}\varphi _0\right\|_H+\\+ \left\| \sum\limits_{j = 1}^\infty  c_jA_0^2\left( {A_0^2  + \gamma _j^2 I} \right)^{ - 1}A_0^2\cos \left( {A_0t} \right)\varphi _0\right\|_H+\\+\left\| \sum\limits_{j = 1}^\infty  c_j\gamma_jA_0^2\left( {A_0^2  + \gamma _j^2 I} \right)^{ - 1}A_0\sin \left( {A_0t} \right)\varphi _0\right\|_H+
\end{multline*}
\begin{multline}\label{eq:mest}
+\left\| \sum\limits_{j = 1}^\infty  c_j\gamma_jA_0^2\left( {A_0^2  + \gamma _j^2 I} \right)^{ - 1}e^{ - \gamma _j t}\varphi _0\right\|_H+\\+ \left\| \sum\limits_{j = 1}^\infty  c_j\gamma_jA_0^2\left( {A_0^2  + \gamma _j^2 I} \right)^{ - 1}\cos \left( {A_0t} \right)\varphi _0\right\|_H+\\+\left\| \sum\limits_{j = 1}^\infty  c_jA_0^2\left( {A_0^2  + \gamma _j^2 I} \right)^{ - 1}A_0\sin \left( {A_0t} \right)\varphi _0\right\|_H\lesssim \left(\sum\limits_{j = 1}^\infty  c_j\right)\left\|A_0^2\varphi _0\right\|_H+\left(\sum\limits_{j = 1}^\infty  c_j\right)\left\|A_0\varphi _1\right\|_H.
\end{multline}
Strictly analogously we obtain the estimate for vector-valued function $Bg_3(t)$
$$
\left\|Bg_3(t)\right\|_H\leqslant \left\|A_0^2g_3(t)\right\|_H\leqslant const \left( \left\|A_0^2\varphi _0\right\|_H+\left\|A_0\varphi_1\right\|_H\right ).
$$
Uniting inequalities \eqref{eq:g1}--\eqref{eq:mest} we have the following estimate
So for the following estimate is valid
\begin{equation}\label{eq:lest}
\left\| {g_1(t)} \right\|_H+\left\| A^2g_2(t)\right\|_H+\left\|Bg_3(t) \right\|_H \leqslant const \left( {\left\| {A_0^2\varphi _0 } \right\|_H  + \left\| {A_0\varphi _1 } \right\|_H } \right).
\end{equation}
From the rpresentation \eqref{eq:h} and estimate \eqref{eq:mest} we obtain that vector-valued functions $A^2g_2(t)$, $Bg_3(t)$ belong to the space $L_{2,\gamma } \left( {\R_+ ,H} \right)$ with arbitraty $\gamma>0$ and the following estimate 
\begin{multline}\label{eq:des}
\left\| {h'(t)} \right\|_{L_{2,\gamma } \left( {\R_+ ,H} \right)}\leqslant \\\leqslant
\left\| {g_1(t)} \right\|_{L_{2,\gamma } \left( {\R_+ ,H} \right)}+\left\|A^2g_2(t)\right\|_{L_{2,\gamma } \left( {\R_+ ,H} \right)}+\left\|Bg_3(t) \right\|_{L_{2,\gamma } \left( {\R_+ ,H} \right)} \leqslant \\\leqslant const \left( {\left\| {A_0^2\varphi _0 } \right\|_H  + \left\| {A_0\varphi _1 } \right\|_H } \right)
\end{multline}
is valid. Let us put $\gamma=\gamma_0$.

Because of the kernels $K(t)$ and $Q(t)$ belong to the space $W_1^1(\mathbb R_+)$, due to Remark \ref{r1} and conditions of the theorem \ref{T:2} vector-valued fucntion $h(t)$ (see \eqref{eq:h}) satisfies the desired property $\mathop {\lim }\limits_{t \to  + 0} h(t) = h( + 0) = 0$. 

At least on the base of estimate \eqref{est13}, repesentations \eqref{eq:h}, \eqref{eq:h'} and estimate \eqref{eq:des} we have the desired estimate \eqref{eq:92}. Theorem \ref{T:2} is proved.

\section{Conclusion remarks.}
In our previous works \cite{16}--\cite{18}, \cite{418} the problem  \eqref{eq:101}, \eqref{eq:102} in the case when operator $B=0$ was investigated in details. The main attention was payed in these articles to the spectral analysis of operator-function \eqref{eq:19} along with the problem of correct solvability. In turn of the base of spectral analysis the results about representation of the storng solutions of the equation \eqref{eq:101} in the series of exponentials corresponding to the points of spectra for operator-valued function $L(\lambda)$ were proved (see \cite{OT} and monograph \cite{418}, chapter 3, for details).

\vspace{0.5ex}

\begin{flushleft}

{\footnotesize N. A. Rautian

Department of Mechanics and Mathematics

Moscow Lomonosov State University

Vorobievi Gori, Moscow, 117234 

Russia}

{\footnotesize { \it Email address:} nrautian@mail.ru}

\end{flushleft}

\vspace{0.5ex}

\begin{flushleft}

{\footnotesize V. V. Vlasov

Department of Mechanics and Mathematics

Moscow Lomonosov State University

Vorobievi Gori, Moscow, 117234 

Russia}

{\footnotesize { \it Email address:} : vlasovvv@mech.math.msu.su}

\end{flushleft}

\end{document}